\documentclass[11pt]{amsart}

\usepackage{amsthm,amsmath,amssymb,amsxtra,amscd}
\usepackage{url}
\usepackage[mathscr]{euscript}
\usepackage[shortlabels]{enumitem}

\newtheorem{theorem}{Theorem}[section]
\newtheorem{lemma}[theorem]{Lemma}
\newtheorem{corollary}[theorem]{Corollary}

\numberwithin{equation}{section}

\renewcommand{\pmod}[1]{\ \left({\rm mod\ } #1 \right)}

\begin{document}

\title{ON COLOURING ORIENTED GRAPHS OF LARGE GIRTH}

\author{P. Mark Kayll}
\address{Department of Mathematical Sciences, University of Montana, 32 Campus Dr, Missoula MT 59812-0864}
\email{mark.kayll@umontana.edu}

\author{Michael Morris}
\address{HelioCampus, 7315 Wisconsin Ave, Suite 750W, Bethesda MD 20814}
\email{michael.morris@alumni.brown.edu}

\date{18 July 2023}



\thanks{\noindent PMK was partially supported by a Simons Foundation grant (\#279367 to Mark Kayll).}
\thanks{This work forms part of MM's MA thesis.}
\thanks{\bfseries To appear in \textit{Contributions to Discrete Mathematics}}

\begin{abstract}
We prove that for every oriented graph $D$ and every choice of
positive integers $k$ and $\ell$, there exists an oriented graph $D^*$
along with a surjective homomorphism $\psi\colon V(D^*) \to V(D)$ such
that: (i) girth$(D^*) \geq\ell$; (ii) for every oriented graph $C$
with at most $k$ vertices, there exists a homomorphism from $D^*$ to
$C$ if and only if there exists a homomorphism from $D$ to $C$; and
(iii) for every $D$-pointed oriented graph $C$ with at most $k$
vertices and for every homomorphism \linebreak $\varphi\colon V(D^*) \to V(C)$
there exists a unique homomorphism $f\colon V(D) \to V(C)$ such that
$\varphi=f \circ \psi$.  Determining the oriented chromatic number of
an oriented graph $D$ is equivalent to finding the smallest integer
$k$ such that $D$ admits a homomorphism to an order-$k$ tournament, so
our main theorem yields results on the girth and oriented chromatic
number of oriented graphs.  While our main proof is probabilistic
(hence nonconstructive), for any given $\ell\geq 3$ and $k\geq 5$, we
include a con\-struction of an oriented graph with girth $\ell$ and
oriented chromatic number $k$.
\end{abstract}

\maketitle

\section{Introduction}

In 1959, Paul Erd\H{o}s~\cite{erdos} famously proved probabilistically
the existence of graphs of arbitrarily large girth and arbitrarily
large chromatic number.  We briefly discuss the history of this and
related topics and point the reader to \cite{Harut} or
\cite{kayllparsa} for more details and references.  We begin by noting
that with \cite{erdos} not suggesting a way to construct such graphs,
one is led to seek such constructions. Other natural directions of
inquiry are (1) generalizing Erd\H{o}s' result and (2) developing
analogues of his results for other types of graphs, specifically of
interest here, directed graphs.

Both refinements and generalizations of \cite{erdos} have followed in
the intervening six-plus decades.  In 1976, Bollob\'{a}s and
Sauer~\cite{bollobas} refined Erd\H{o}s' result by showing that for
any positive integer $n$ there are graphs of arbitrarily large girth
that are `uniquely' $n$-colourable.  In 1996, Zhu~\cite{zhu}, working
with graph homomorphisms as a generalization of colouring, was able to
carry forward the work of \cite{bollobas} by showing that for any
`core' $H$, there are uniquely $H$-colourable graphs of arbitrarily
large girth.  As complete graphs are cores, Zhu's work generalizes
\cite{erdos}.  Zhu's main result in \cite{zhu} was further generalized
by Ne\v{s}et\v{r}il and Zhu~\cite{nestril} to the notion of `pointed'
graphs.  We follow a similar trajectory in the present paper. 

So let us shift our attention to digraphs. Bokal et al.~\cite{circ}
studied the digraph circular chromatic number and showed that digraph
colouring theory is similar to that of undirected graphs.  For
undirected graphs $G$, the circular chromatic number $\chi_c(G)$ is a
refinement of the chromatic number $\chi(G)$ because
$\chi(G)-1 < \chi_c(G) \leq \chi(G)$ (see, e.g., \cite{circz}). 
Analogously, the circular chromatic number $\chi_c(D)$ of a digraph
refines the chromatic number $\chi(D)$, here defined to be the minimum
integer $k$ such that $V(D)$ can be partitioned into $k$ acyclic
subsets. The former parameter is defined using using `acyclic'
homomorphisms---see \cite{circ} for details---which 
introduced complications.  For example, the authors of
\cite{kayllparsa} had to use a lot of care to demonstrate that certain
mappings don't fail to be acyclic homomorphisms.  The fact that we
consider oriented colouring here means we have no need to turn to
acyclic homomorphisms.  The reader might appreciate how much this
simplifies our proofs in comparison with those of \cite{kayllparsa}. 

Subsequently to \cite{circ}, a subset of the authors and their
doctoral students in \cite{Harut} completed work in the realm of
digraphs analogous to that of Zhu for graphs in \cite{zhu}.  Then
\cite{kayllparsa} generalized the results of \cite{circ,Harut} just as
Ne\v{s}et\v{r}il and Zhu in \cite{nestril} generalized
\cite{erdos,zhu}.  One of our successes in the present work is a
similar sequence of generalizations for oriented graphs. 

We delay definitions for a little longer (until Section~\ref{TandN})
and proceed to state our main result and a couple of its consequences: 

\begin{theorem}
\label{mainresult} 
For every oriented graph $D$ and every choice of positive integers $k$
and $\ell$, there exists an oriented graph $D^*$ along with a
surjective homomorphism $\psi\colon  V(D^*) \to V(D)$ such that: 
\begin{enumerate}[(i)]
    \item girth$(D^*) \geq \ell$;
    \item for every oriented graph $C$ with at most $k$ vertices,
      there exists a homo\-morph\-ism from $D^*$ to $C$ if and only if
      there exists a homo\-morphism from $D$ to $C$; and 
    \item for every $D$-pointed oriented graph $C$ with at most $k$
      vertices and for every homomorphism
      $\varphi\colon  V(D^*) \to V(C)$ there exists a unique
      homomorphism $f\colon  V(D) \to V(C)$ such that $\varphi=f \circ \psi$. 
\end{enumerate}
\end{theorem}

An attentive reader familiar with \cite{kayllparsa} may be concerned
that our Theor\-em~\ref{mainresult} is an immediate consequence of
\cite[Theorem~1]{kayllparsa}.  After all, the earlier result applies
to digraphs in general, and oriented graphs are a specific type of
digraph.  Furthermore, oriented colourings are homomorphisms from
oriented graphs to oriented graphs, so in particular they are acyclic
homo\-morphisms.  We can see this because preimages under a
homomorphism must be independent sets, and are hence acyclic.  The
important difference is that in our Theorem~\ref{mainresult}, we are
able to get an \emph{oriented} $D^*$ and an \emph{oriented colouring}
$\psi$, whereas \cite{kayllparsa} guarantees only a \emph{digraph}
$D^*$ and an \emph{acyclic homomorphism} $\psi$.  Although $D^*$ of
\cite{kayllparsa} will in fact be an oriented graph when $D$ is an
oriented graph, one can readily check that the acyclic homomorphism
$\psi$ of \cite{kayllparsa} in general will not be an oriented
colouring, so the earlier results do not guarantee the desired results
for oriented graphs. The importance of this distinction becomes clear
as we discuss two consequences of Theorem~\ref{mainresult}, which we
now state.

\begin{corollary}
\label{Erdos} If $D$ and $C$ are oriented graphs such that $D$ is not
$C$-colourable, then for every positive integer $\ell$, there exists an
oriented graph $D^*$ of girth at least $\ell$ that is $D$-colourable
but not $C$-colourable. 
\end{corollary}

\begin{corollary}
\label{Zhu}For every oriented core $D$ and every positive integer
$\ell$, there is an oriented graph $D^*$ of girth at least $\ell$ that
is uniquely $D$-colourable. 
\end{corollary}
 
To see that Theorem~\ref{mainresult} implies Corollary~\ref{Erdos}, if
we have $D$ and $C$ as in Corollary~\ref{Erdos} with a given integer
$\ell$ and take $k$ to be the order of $C$, then (i) of
Theorem~\ref{mainresult} gives us a $D^*$ of required girth such that
$\psi\colon  V(D^*) \to V(D)$, so $D^*$ is $D$-colourable.  But as $D$
is not $C$-colourable, condition (ii) of Theorem~\ref{mainresult}
implies that $D^*$ is not $C$-colourable. 

To see that Theorem~\ref{mainresult} implies Corollary~\ref{Zhu}
follows a similar argument as in \cite{kayllparsa}.  We note that
cores $D$ are $D$-pointed.  So if we are given a positive integer
$\ell$ and a core $D$, we can take $k=|V(D)|$.  Then
Theorem~\ref{mainresult} gives a $D^*$ of girth at least $\ell$ and a
$D$-colouring $\psi\colon  V(D^*) \to V(D)$.  We can set $C=D$ in part
(iii) of Theorem~\ref{mainresult}, which gives us that for every
$D$-colouring $\varphi\colon  V(D^*) \to V(D)$ there is a (unique)
homomorphism $f\colon  V(D) \to V(D)$ such that
$\varphi=f \circ\psi$.  Because $D$ is a core, $f$ is an automorphism,
so $\varphi$ and $\psi$ differ by this automorphism and $D^*$ is
indeed uniquely $D$-colourable. 

\section{Terminology and notation}
\label{TandN}

We assume basic familiarity with graphs and digraphs and refer the
reader to \cite{bondy} for any omitted concepts.
Here we  consider oriented graphs and oriented colourings going forward
unless indicated otherwise.  An \emph{oriented graph} $D$ is a digraph
in which for every pair of vertices $u,v$, at most one of $uv$ and
$vu$ is an element of $A(D)$, the \emph{arc set} of $D$.  Our oriented
graphs will always be finite and loopless without multiple arcs;
opposite arcs are precluded by the definition of oriented graphs.  It
can be easier to think about an oriented graph as one obtained by
assigning directions to each edge of some (undirected) graph $G$.
Recall that a \emph{tournament} $D$ on $n$ vertices is an oriented
graph obtained by assigning a direction to each edge of the complete
graph $K_n$. \emph{Cycles} of oriented graphs are directed cycles, and
the \emph{girth} of an oriented graph $D$ is the length of a shortest
directed cycle in $D$.  Finally, for oriented graphs $D$ and $C$, an
\emph{oriented graph homomorphism} is a map $f\colon  V(D) \to V(C)$
such that whenever $xy \in A(D)$, we also have $f(x)f(y)\in A(C)$. 

We are now ready to define an `oriented colouring' of an oriented graph
$D$. An \emph{oriented $k$-colouring}, then, is a map 
$c\colon  V(D) \to \{1,\ldots,k\}$ such that:
\begin{enumerate}
    \item $c(x)\neq c(y)$ for every arc $xy \in A(D)$, and
    \item $c(u) \neq c(y)$ for every two arcs $uv \in A(D)$ and
      $xy \in A(D)$ with $c(v)=c(x)$. 
\end{enumerate}
This is by now a standard definition; see, e.g., \cite{oriented}.

This definition of an oriented colouring is equivalent to that of a
homo\-morphism to a tournament on $k$ vertices. First, it is clear that
a homomorph\-ism to a tournament satisfies condition (1) of being an
oriented colouring because it is a homomorphism, and condition (2) is
satisfied because tour\-naments have no opposite arcs.  On the other
hand, given such a map $c$, we can construct an oriented graph $C^*$
with $V(C^*)=\{1,\ldots,k\}$ and
$A(C^*)=\{xy\colon x,y \in V(C^*)\text{ and } xy=c(a)c(b) \text{ for some }ab \in A(D)\}$.
Then it is clear that $C^*$ is a subgraph of a tournament $C'$ on $k$
vertices by property (2) of $c$.  Furthermore, $C^*$ was constructed
so that $c$ is a homomorphism to $C^*$ and thus a homomorphism to
$C'$, so $c$ is a homomorphism to a tournament on $k$ vertices.  We
always consider oriented colourings to be homomorphisms to tournaments.   

Having defined an oriented colouring, we now give the related
definition of `oriented chromatic number.'  Given an oriented graph
$D$, its \emph{oriented chromatic number} $\chi_{o}(D)$ is the minimum
number of vertices of an oriented graph $C$ such that there exists a
homomorphism of $D$ to $C$.  As $C$ is always a subgraph of some
tournament $T$, we will always consider the oriented chromatic number
of $D$ as the minimum number of vertices of a tournament $T$ such that
there exists a homomorphism of $D$ to $T$. 

For terminology more directly related to our theorem statements, we
say that a homomorphism of oriented graphs of $D$ to $C$ is a
\emph{$C$-colouring of} $D$, and we say that $D$ is
\emph{$C$-colourable}.  We say that $D$ is \emph{uniquely}
\emph{$C$-colourable} if there is a homomorphism of $D$ onto $C$, and
for any two $C$-colourings $\psi$ and $\varphi$ of $D$, these
homomorphisms `differ by an automorphism'.  That is, there is some
$f\in \text{Aut}(C)$ such that $\psi=f \circ \varphi$.  For an oriented
graph $D$, we say that $D$ is a \emph{core} if every homomorphism
$f\colon  V(D) \to V(D)$ is an automorphism.  Finally, we say that for
oriented graphs $C$ and $D$, the digraph $C$ is $D$\emph{-pointed} if
there do not exist two distinct $C$-colourings of $D$ that agree on all
but one vertex of $D$. 

\section{ Setup for the proof of Theorem~\ref{mainresult}}

For a given oriented graph $D$, we begin the `construction' of the
digraph $D^*$, and we do so by first constructing a digraph $D_0$,
again inspired by \cite{kayllparsa}. 

We define $V(D_0)=V_1 \cup V_2 \cup \cdots \cup V_a$ where
$V(D)=\{1,2,\ldots,a\}$, and each $|V_i|=n$ for some fixed $n$ large
enough to satisfy necessary probabilistic inequalities.  Then we
define the arc set
$A(D_0)=\{xy: x \in V_i, y \in V_j \text{ and } ij\in A(D)\}$.  We can
view each $V_i$ simply as the preimage of a vertex $i \in V(D)$ under
the natural homomorphism $\psi\colon V(D_0) \to V(D)$, mapping each
$V_i$ to $i$, for $i \in \{1,\ldots,a\}$.   

Now we use $D_0$ to `construct' an oriented graph $D^*$
probabilistically.  First we fix an $\varepsilon$ with
$0<\varepsilon<1/(4\ell)$ where $\ell$ is chosen as in the statement
of Theorem~\ref{mainresult}.  Then our random oriented graph model
$\mathscr{D}(n,p)$ 
consists of spanning subgraphs of $D_0$ where arcs are chosen randomly
and independent\-ly with probability $p=n^{\varepsilon-1}$ with $n$
sufficiently large.  We now introduce three lemmas from
\cite{kayllparsa}. 

\begin{lemma}
\label{cycles}
(i) The expected number of cycles of length less than $\ell$ in a
digraph
$\hat{D} \in \mathscr{D}(n,p)$ is bounded above by
$n^{\varepsilon \ell}n^{-\varepsilon/2}$;\\
(ii) the expected number of pairs of cycles of length less than $\ell$
in a digraph $\hat{D} \in \mathscr{D}(n,p)$ which intersect in at
least one vertex is bounded above by $n^{-1/2}$. 
\end{lemma}

This is Lemma~5 of \cite{kayllparsa}, except that our oriented graph
model $\mathscr{D}(n,p)$ differs.  In particular, our $D_0$ has fewer
arcs than the analogue in \cite{kayllparsa}, so the lemma remains true
in our case. This along with Markov's Inequality shows that
asymptotically almost all oriented graphs in $\mathscr{D}(n,p)$ have
at most $n^{\varepsilon \ell}$ cycles of length less than $\ell$ which
are pairwise vertex-disjoint; see, e.g., \cite{kayllparsa}. 

We introduce some definitions from \cite{kayllparsa} (which itself
adopted these from \cite{nestril}), first calling a set
$\mathscr{A} \subseteq V(D_0)$ \emph{large} if there are distinct
$i,j \in\{1,\ldots,a\}$ with $ij \in A(D)$ such that
$|\mathscr{A} \cap V_i|\geq n/k$ and $|\mathscr{A} \cap V_j| \geq n/k$,
and calling $ij \in A(D)$ in this case a \emph{good arc} for
$\mathscr{A}$.  Then given a large $\mathscr{A}$, we denote by
$|\hat{D}/\mathscr{A}|$ the minimum number of arcs of a random
$\hat{D}$ which lie in the set
$\{xy\colon x \in \mathscr{A} \cap V_i,y \in \mathscr{A} \cap V_j\}$,   
taken over all instances in which $ij$ a good arc. Then we have: 

\begin{lemma}[\cite{kayllparsa}]
  \label{large}
If $\hat{D} \in \mathscr{D}(n,p)$ and $\mathscr{A}$ is
large, then $P(|\hat{D}/\mathscr{A}| \geq n)=1-o(1)$. 
\end{lemma}

Again the space $\mathscr{D}(n,p)$ in \cite{kayllparsa} differs from
ours, but the proof still follows through unchanged because the arcs
counted in $|\hat{D}/\mathscr{A}|$ in \cite{kayllparsa} are all
present in the current model. 

We shall need to adopt one last lemma from 
\cite{kayllparsa}, and its validity here follows using similar 
arguments to those for Lemmas~\ref{cycles} and \ref{large}.

\begin{lemma}[\cite{kayllparsa}]
\label{arcs} 
For almost all digraphs in $\mathscr{D}(n,p)$, 
all nonempty $\mathscr{A} \subseteq V_{v}$ and 
$\mathscr{B} \subseteq V_{i_0}$ (for $v,i_0 \in\{1,\ldots,a\}$
with $vi_0 \in A(D)$) with 
$|\mathscr{A}|=n-|\mathscr{B}|(k-1)$ and 
$|\mathscr{B}|\leq n/k$ satisfy the property of 
$\mathscr{A} \cup \mathscr{B}$ inducing more than $\text{min}\{|\mathscr{B}|,n^{\epsilon \ell}\}$ 
arcs from $\mathscr{A}$ to $\mathscr{B}$.
\end{lemma}

\noindent
Now we can move on to the proof of our main theorem.

\section{Proof of Theorem~\ref{mainresult}}

Lemma~\ref{cycles} and its consequences mean that asymptotically
almost all $D'\in \mathscr{D}(n,p)$ have at most
$n^{\varepsilon\ell}$ pairwise-disjoint cycles of length less than $\ell$.
Similarly, Lemma~\ref{large} guarantees that asymptotically almost all
$D'\in \mathscr{D}(n,p)$ have the property that all good arcs of $D$
for large sets $\mathscr{A}$ induce at least $n$ arcs of $D'$.
Finally Lemma~\ref{arcs} guarantees the existence of necessary arcs as
described later for almost all $D'$.  Therefore, there exists some
$D'\in \mathscr{D}(n,p)$ enjoying the three stated properties, and we
select such a $D'$.  Now we pick one arc from each of the at most
$n^{\varepsilon\ell}$ cycles of length less than $\ell$ in $D'$,
giving an independent arc set (i.e., a matching)  
$M$, and define $D^*=D'-M=(V(D_0),A(D')\setminus M)$.  
It is clear then that $D^*$ has girth at least $\ell$, 
and that $\psi\colon V(D^*) \to V(D)$ defined by 
$\psi(x)=i$ if and only if $x \in V_i$ gives a surjective 
homomorphism, yielding (i) from Theorem~\ref{mainresult}.  Note that
since $\varepsilon<1/(4\ell)$, the deleted arc set 
satisfies $|M|\leq n^{\varepsilon \ell} < n^{1/4}$. 

Now we work toward (ii) from Theorem~\ref{mainresult}.  Let us fix an
oriented graph $C$ of order at most $k$, and assume that there is a
homomorphism $\varphi\colon  V(D^*) \to V(C)$.  Then the Pigeonhole
Principle implies that for every $i \in V(D)$ there is a vertex
$x \in V(C)$ such that $|V_i \cap \varphi^{-1}(x)| \geq n/k$.  Then let us
define $f\colon  V(D) \to V(C)$ by $f(i)=x$ for some $x \in V(C)$ such
that $|V_i \cap \varphi^{-1}(x)| \geq n/k$.  We must show this $f$ is
a homomorphism. 

Let $ij \in A(D)$ and consider all possible $a,b \in V(D^*)$ where
$a\in V_i \,\cap\, \varphi^{-1}(f(i))$ and $b\in V_j\cap\varphi^{-1}(f(j))$.
If there is one such arc $ab \in A(D^*)$, this
will guarantee the existence of an arc $f(i)f(j) \in A(C)$ by the
existence of $\varphi$.  Recall that $f$ satisfies
$|V_i \cap \varphi^{-1}(f(i))| \geq n/k$ and
$|V_j \cap \varphi^{-1}(f(j))| \geq n/k$.  Then
$\mathscr{A}=\left(V_i \cap \varphi^{-1}(f(i))\right)\cup\left(V_j\cap \varphi^{-1}(f(j))\right)$
is large as defined for
Lemma~\ref{large}, so by our choice of $D'$ relying on that lemma,
$D'$ has at least $n$ arcs with endpoints in $\mathscr{A}$.  Then
since we have removed at most $n^{1/4}$ arcs from $D'$ to construct
$D^*$, there exists at least one such arc $ab \in A(D^*)$, and in fact
many such arcs.  So we have $\varphi(a)\varphi(b) \in A(C)$, and we
have that $f(i)=\varphi(a)$ and $f(j)=\varphi(b)$ with
$f(i) \neq f(j)$ because $\varphi$ is a homomorphism.  So $f(i)f(j) \in A(C)$,
and $f$ maps arcs to arcs and is thus a homomorphism. 

Conversely, if we assume that there is a homomorphism
$f\colon  V(D) \to V(C)$, then we get the homomorphism
$\varphi\colon V(D^*) \to V(C)$ by $\varphi=f \circ \psi$,
completing our proof of (ii). 

Now we turn to (iii), letting $C$ be a $D$-pointed oriented graph of
order at most $k$, and $\varphi\colon V(D^*) \to V(C)$ be a
homomorphism.  We shall use $f\colon  V(D) \to V(C)$ as in the proof
of (ii).   The $D$-pointedness of $C$ forces for every $i \in V(D)$
the existence of a unique $x_i \in V(C)$ such that $|\varphi^{-1}(x_i)
\cap V_i| \geq n/k$. We demonstrate this using an argument similar to
that in \cite{kayllparsa}. If some $x_i$ were not unique and $x_i'$
also satisfies $|\varphi^{-1}(x_i') \cap V_i| \geq n/k$, then we could
define $f'$ by  
\[
f'(j)=\left\{\begin{array}{cl}
f(j) & \text{for } j \neq i\\
x_i' &\text{for } j = i,
\end{array}\right.
\]
giving another homomorphism differing at one vertex of $D$ and
contradicting the $D$-pointed\-ness of $C$.  This establishes the
uniqueness of a homomorphism $f$ chosen in this way.  If we assume
that $\varphi \neq f \circ \psi$, then there must be some vertex
$z\in V(D^*)$ such that $\varphi(z) \neq (f \circ \psi)(z)$.  So if
$z\in V_j$, then  
$\varphi(z)\neq (f \circ \psi)(z)=f(j)$.  
Thus $V_j\setminus\left(\varphi^{-1}(f(j))\cap V_j\right)\neq\varnothing$
(as it contains $z$), which leads to a contradiction as we proceed to show.

We begin by choosing a vertex $i_0 \in \{1,\ldots,a\}$ so that
$t:=|\varphi^{-1}(f(i_0))\cap V_{i_0}|$ is minimized; the definition
of $f$ gives $t\geq n/k$ while the purported 
$z$ of the preceding paragraph gives $t<n$. The last inequality
shows that $\varphi^{-1}(f(i_0))\cap V_{i_0}$
is a proper subset of $V_{i_0}$. Let us now choose a vertex
$x\in V(C)$, distinct from $f(i_0)$, so as to
maximize the size of the set 
$\mathscr{B}:=\varphi^{-1}(x) \cap V_{i_0}$. Denoting
this size by $b=|\mathscr{B}|$, we see that $b<n/k$ by the 
previously established uniqueness property of $f$
(exactly one vertex of $V(C)$ satisfies `$\geq n/k$' here, and
$f(i_0)\neq x$ is already that witness). 
Notice also that these new parameters satisfy
\begin{equation}
    \label{b-t-relation}
    b(k-1)\geq n-t 
\end{equation}
because the (at most) $(k-1)$
preimages $\varphi^{-1}(y)\cap V_{i_0}$ within $V_{i_0}$
(as $y$ runs through $V(C)\setminus\{f(i_0)\}$)
exhaust the $(n-t)$ vertices within $V_{i_0}$
that are not mapped to $f(i_0)$ by $\varphi$.

Now we define $f'\colon V(D) \to V(C)$ by:
\[
f'(i)=\left\{\begin{array}{cl}
f(i) & \text{for } i \neq i_0\\
x &\text{for } i = i_0.
\end{array}\right.
\]
\noindent
Because $f$ and $f'$ differ only at $i_0$ and $C$ is
$D$-pointed, $f'$ is not a homomorph\-ism.  Thus, it fails 
to send arcs to arcs.  So it must be for some $v \in V(D)$,
distinct from $i_0$, either $vi_0 \in A(D)$ and 
$f(v)x \notin A(C)$ or $i_0v \in A(D)$ and 
$xf(v) \notin A(C)$.  Without loss of generality, we 
assume that $vi_0\in A(D)$ and 
\begin{equation}
    \label{for-contradiction}
    f(v)x \notin A(C).
\end{equation}
With $v$ being among the candidate vertices $1,\ldots,a$ 
during our choice of $i_0$, we have 
$|\varphi^{-1}(f(v))\cap V_{v}|\geq t$, and 
(\ref{b-t-relation}) shows that $t\geq n-b(k-1)$;
therefore, we can select a subset
$\mathscr{A}\subseteq\varphi^{-1}(f(v))\cap V_v$ with
$|\mathscr{A}|=n-b(k-1)$. 

Because we chose a digraph $D'$ 
satisfying the likely properties articulated in 
Lemma~\ref{arcs}, we know that there are more than 
$\min\{|\mathscr{B}|,n^{\epsilon \ell}\}$ arcs from 
$\mathscr{A}$ to $\mathscr{B}$ in $D'$.
And because the arcs removed from $D'$ to form $D^*$  
comprised a matching 
of size at most $n^{\epsilon \ell}$,
no matter which entry achieves 
$\min\{|\mathscr{B}|,n^{\epsilon \ell}\}$,
there exist vertices $a \in \mathscr{A}$ and 
$b \in \mathscr{B}$ such that $ab \in A(D^*)$. 
Lastly, because $\varphi$ is a homomorphism, we have
$\varphi(a)\varphi(b) \in A(C)$.  But $\varphi(a)=f(v)$ 
and $\varphi(b)=x$, so that $f(v)x \in A(C)$, which
contradicts (\ref{for-contradiction}).  Therefore, our assumption of the 
existence of a vertex $z \in V(D^*)$ such that 
$\varphi(z)\neq (f\circ\psi)(z)$ is incorrect, 
and we must have $\varphi=f\circ\psi$.  Finally, 
we note that the surjectivity of $\psi$ implies that such a
homomorphism $f$ is unique. \qed

\section{Constructions}

Our last natural direction of exploration from Erd\H{o}s' original
theorem is that of actually constructing those graphs which we have
probabilistically proven exist.  These constructions are generally
challenging and delicate.  The common approach is to proceed by
induction, constructing a (di)graph of chromatic number $n+1$ with
girth $\ell$ using copies of a (di)graph of chromatic number $n$ with
girth $\ell$.  The first such construction was completed by
Lov\'asz~\cite{lovasz} in 1968 using hypergraphs intermediately.  It
was not until 1989 that K\v{r}\`{i}\v{z}~\cite{kriz} was able to
create a purely graph-theoretic construction of highly chromatic
graphs without short cycles.  Similarly, Severino~\cite{cons1}
demonstrated constructions of highly chromatic digraphs without short
cy\-cles and in \cite{cons2} constructed uniquely $n$-colourable digraphs
with arbitrarily large girth.

Ideally, we would like to construct the digraph $D^*$ with all the
properties described in Theorem~\ref{mainresult}.  We shall content
ourselves with a construction of an oriented graph of a given girth
and oriented chromatic number and leave the construction of such a
$D^*$ for future authors.

\begin{theorem}
\label{const}
For integers $k \geq 5$ and $\ell \geq 3$, there exists an oriented
graph $D$ with oriented chromatic number $k$ and girth $\ell$.
\end{theorem}

\noindent
\emph{Remark}: Some instances of $(k,\ell)$ with $k=3$ or $k=4$ are
also feasible.  In particular, $k=3$ is feasible for
$\ell \equiv 0\pmod{ 3}$, and $k=4$ is feasible for all $\ell \geq 3$
with $\ell \neq 5$.  However, we state Theorem~\ref{const} as such
because when $\ell=5$, of necessity our basis starts at $k=5$.
Readers may find it illustrative to convince themselves that the
directed $5$-cycle admits no homomorphism to a tournament on four
vertices, while a directed cycle of any other order admits such a
homomorphism.

\begin{proof} 
We follow the common approach to which we alluded above and proceed by
induction on $k$, so let us fix integers $k$ and $\ell$.  Then we
begin by considering $\overrightarrow{C_\ell}$, an oriented cycle of
length $\ell$ (and girth $\ell$).  We define
$V(\overrightarrow{C_\ell})=\{v_0,v_1,\ldots,v_{\ell-1}\}$, and there
is a homomorphism $c\colon  V(\overrightarrow{C_\ell}) \to V(T_5)$
where $V(T_5)=\{t_0,t_1,t_2,t_3,t_4\}$ and
$\{t_0t_1,t_1t_2,t_2t_3,t_3t_4,t_2t_0,t_3t_0,t_4t_0\} \subseteq A(T_5)$.
Then if $\ell \equiv 0 \pmod{3}$, we have
$c\colon V(\overrightarrow{C_\ell}) \to V(T_5)$ defined by 
\[c(v_r)=t_{r \bmod{3}}.\]
If $\ell \equiv 1 \pmod{3}$, then $c$ is defined by
\[
c(v_r)=\left\{\begin{array}{cl}
t_{r \bmod{3}} & \text{for } r<\ell-1\\
t_3 &\text{for } r=\ell-1.
\end{array}
\right.
\]
And finally, if $\ell \equiv 2 \pmod{3}$, then $c$ is defined by
\[
c(v_r)=\left\{\begin{array}{cl}
t_{r \bmod{3}}& \text{for } r<\ell-2\\
t_3 &\text{for } r=\ell-2\\
t_4 & \text{for } r= \ell-1.
\end{array}
\right.
\]
We note that our base cases have given us an oriented graph of girth
$\ell$ with oriented chromatic number $k \leq 5$.  The verification of
our induction below will then guarantee the existence of an oriented
graph of girth $\ell$ with any given oriented chromatic number
$k \geq 5$.

Having established our base cases, we now proceed with the induction.
So assume we have an oriented graph $D_k$ of girth $\ell$, 
oriented chromatic number $k$, and order $m$, and then define
$V\left(D_k\right)= \{v_0,v_1,\ldots,v_{m-1}\}$.  Because $D_k$ has 
oriented chromatic number $k$, there exists a tournament $T_k$ with
$V(T_k)=\{t_0,t_1,\ldots,t_{k-1}\}$ and a homomorphism
$\varphi_k\colon  V(D_k) \to V(T_k)$.  Now we construct $D_{k+1}$
and the corresponding $T_{k+1}$.  Define the vertex set
$V(D_{k+1})=V(D_k)\cup \{v_{m}\}$, and define the arc set 
$A(D_{k+1})=A(D_k)\cup \{v_iv_{m}\colon i \in \{0,1,\ldots,m-1\}\}$.  
Then we construct $T_{k+1}$ in exactly the same fashion; i.e.,
$V(T_{k+1})=V(T_k)\cup \{t_{k}\}$ and 
$A(T_{k+1})=A(T_k)\cup \{t_it_{k}\colon i \in \{0,1,\ldots,k-1\}\}$.

We now examine the girth and oriented chromatic number of $D_{k+1}$.
First, it is immediately clear that we have created no new oriented
cycles in this construction, so $D_{k+1}$ also has girth $\ell$.  It
is equally clear that we have a homomorphism
$\varphi_{k+1}\colon V(D_{k+1}) \to V(T_{k+1})$ defined by 
\[
\varphi_{k+1}(v)=\left\{\begin{array}{cl}
\varphi_k(v) & \text{for } v \neq v_{m}\\
t_{k} & \text{for } v=v_{m}.
\end{array}
\right.
\]
Therefore, $\chi_o(D_{k+1}) \leq k+1$.

To complete the proof, it remains to show that $D_{k+1}$ admits no
homo\-morphism to a tournament on $k$ vertices.  Assume to the
contrary that for some order-$k$ tournament $T'_k$ the digraph
$D_{k+1}$ admits a homomorphism $\psi\colon  V(D_{k+1}) \to V(T'_k)$.
Let's say that $\psi(v_{m})=x \in V(T'_k)$.  Then because every vertex
$v \in V(D_{k+1})\setminus \{v_m\}$ forms an arc $vv_{m}$, we know
that $\psi(v)\neq x$ for every $v \neq v_{m}$.  If we let $\Lambda$ be
the subgraph of $D_{k+1}$ induced by the vertex set
$\{v_0,\ldots,v_{m-1}\}$, then $\Lambda$ is isomorphic to $D_k$.
Similarly, if we let $\Gamma$ be the subgraph of $T'_k$ induced by
$V(T_k')\setminus \{x\}$, then $\Gamma$ is a tournament on $k-1$
vertices.  But then $\psi|_{V(\Lambda)}$ gives a homomorphism from
$\Lambda$ to $\Gamma$, a tournament on $k-1$ vertices, contradicting
the fact that $D_k$ has oriented chromatic number $k$.  Therefore,
$D_{k+1}$ indeed has oriented chromatic number $k+1$.
\end{proof}

\section*{Acknowledgements}
The authors are grateful to Gary MacGillivray for suggesting this
research direction following  a preliminary presentation of
\cite{kayllparsa} at the 2019 Coast Com\-binatorics Conference. Thanks
also to an anonymous referee for the careful reading and, in
particular, for catching a  gap in our original presentation of  
Theorem~\ref{mainresult}'s proof.

\providecommand{\bysame}{\leavevmode\hbox to3em{\hrulefill}\thinspace}
\providecommand{\MR}{\relax\ifhmode\unskip\space\fi MR }
\providecommand{\MRhref}[2]{%
  \href{http://www.ams.org/mathscinet-getitem?mr=#1}{#2}
}
\providecommand{\href}[2]{#2}

\end{document}